\documentclass[reqno, 11pt]{amsart}
\usepackage{amssymb}
\usepackage{amsthm}
\usepackage{amsmath}
\usepackage{times}

\usepackage{latexsym}
\usepackage[mathscr]{eucal}
\usepackage{float}
\numberwithin{equation}{section}

\newtheorem{theorem}{Theorem}[section]
\newtheorem{proposition}[theorem]{Proposition}
\newtheorem{lemma}[theorem]{Lemma}
\newtheorem{corollary}[theorem]{Corollary}

\newtheorem{remark}[theorem]{Remark}

\newcommand{\mainthmname}{}

\newtheorem{definition}[theorem]{Definition}
\newtheorem{example}[theorem]{Example}

\newcommand{\KN}{\mathbin{\bigcirc\mspace{-15mu}\wedge\mspace{3mu}}}

\title[Conformally Sasakian manifolds and the $k$-nullity condition]{Conformally Sasakian manifolds and the $k$-nullity condition}
\author[Fortun\'e Massamba, Pontsho Moile]{Fortun\'{e} Massamba*, Pontsho Moile**}
\newcommand{\acr}{\newline\indent}
\address{\llap{*\,} Discipline of Mathematics\acr 
	School of Agriculture and Science\acr
	University of KwaZulu-Natal\acr
	Private Bag X01, Scottsville 3209\acr 
	South Africa} 
\email{massfort@yahoo.fr, Massamba@ukzn.ac.za} 
\address{\llap{**\,}Discipline of Mathematics\acr 
	School of Agriculture and Science\acr
	University of KwaZulu-Natal\acr
	Private Bag X01, Scottsville 3209\acr 
	South Africa \acr 
	and\acr
	Department of Mathematics\acr
	University of Botswana\acr
	Private Bag 0022\acr
	Gaborone, Botswana}
\email{Moilep@ub.ac.bw} 
\thanks{}
\thanks{}


\begin{document}

\begin{abstract}

We study almost contact metric manifolds that are locally conformal to Sasakian manifolds in the normalization of Olszak. We show that the local conformal factors are uniquely determined by the structure; consequently, every such manifold is globally conformal to a Sasakian manifold and its Lee form is exact. We then investigate the $k$-nullity condition
$
R(X,Y)\xi=k\{\eta(Y)X-\eta(X)Y\}.
$
We prove that a smooth nullity function is necessarily constant, including in dimension three. Moreover, the $k$-nullity condition is equivalent to the assertion that $f=\exp(\sigma)$ satisfies
$
\mathrm{Hess}_g f=-(kf+C)g
$
for a constant $C$, and hence to the fact that the conformal metrics $g$ and $g_0=\exp(-2\sigma)g$ have the same geodesic circles. We derive the range of the conformal factor, a local warped-product description, and a classification in the complete case: when the Lee form is nonzero, the complete metric is a simply connected space form of curvature $k$. Finally, we construct incomplete examples of nonconstant sectional curvature in every dimension $2n+1\geq 5$.	
\end{abstract}
\maketitle

\section{Introduction}\label{sec:intro}

Sasakian manifolds are the odd-dimensional counterparts of K\"ahler manifolds \cite{bl,boyer-galicki}. Their defining identity
$
(\nabla_X\phi)Y=g(X,Y)\xi-\eta(Y)X
$
fixes the curvature along the characteristic vector field,
$
R(X,Y)\xi=\eta(Y)X-\eta(X)Y.
$
It is therefore natural to ask how much of this rigidity survives under conformal changes of the structure. In even dimensions, such changes lead to locally conformal K\"ahler geometry, where the local conformal factors agree only up to constants and the Lee form is closed but not necessarily exact \cite{dragomir-ornea}. In odd dimensions, conformal changes of almost contact metric structures were introduced by Vaisman \cite{vaisman} and studied by Chinea and Marrero \cite{chinea-marrero}, within the classification of Chinea and Gonz\'alez \cite{chinea-gonzalez}, and by Kirichenko and Uskorev \cite{kirichenko-uskorev}. Olszak \cite{olszak} introduced locally conformal almost cosymplectic manifolds, whose curvature under a $(k,\mu)'$-nullity condition was studied by Maduna and Massamba \cite{madu}, see also \cite{mass1}. For a Sasakian local model, Falcitelli \cite{falcitelli1} showed that manifolds locally conformal to $\mathcal{C}_6$-manifolds are, in dimension $2n+1\geq 5$, either locally conformal cosymplectic or globally conformal to Sasakian manifolds. She also studied the $k$-nullity condition in this class. Such a manifold is a generalized Sasakian-space-form in the sense of Alegre, Blair and Carriazo \cite{abc} if and only if it satisfies the $k$-nullity condition and has pointwise constant $\phi$-sectional curvature. Local classifications are also obtained under curvature assumptions. Conformal Sasakian manifolds were also introduced by Abedi \textit{et al.} \cite{abt}, who studied submanifolds of conformal Sasakian space forms, while Abedi and Ilmakchi \cite{abedi} studied submanifolds of conformal Sasakian manifolds.

Nullity conditions measure how far a structure is from the Sasakian model. Tanno \cite{tanno88} introduced the $k$-nullity condition
$
R(X,Y)\xi=k\{\eta(Y)X-\eta(X)Y\}.
$
Blair, Koufogiorgos and Papantoniou \cite{bkp} introduced $(\kappa,\mu)$-spaces, which were classified by Boeckx \cite{boeckx}, and Dacko and Olszak \cite{dacko-olszak} studied their almost cosymplectic analogues. These conditions are flexible. They admit many complete examples that are not of constant curvature, and in dimension $3$ the nullity functions need not be constant \cite{koufogiorgos-tsichlias}. In the almost cosymplectic setting, Massamba \cite{mass2} studied a generalized $(k,\mu)'$-nullity condition, in which $k$ and $\mu$ are smooth functions. He proved that if $h'=\phi\circ h\neq 0$, then $k<0$ and
$
\mu=\xi(\ln\sqrt{-k}),
$
and that if, in addition, $\nabla_{\xi}h'=0$, the manifold is locally isometric to a warped product
$
\mathbb{H}^{n+1}(k)\times_f\mathbb{R}^{n}
$
of a hyperbolic space and a Euclidean space.

In this paper, we study almost contact metric manifolds that are locally conformal to Sasakian manifolds, with Olszak's normalization, and the $k$-nullity condition on them. We first show that the theory is not genuinely local. The local conformal factors are determined pointwise by the structure, so every locally conformal (l.c.) Sasakian manifold is globally conformal to a Sasakian manifold and has an exact Lee form, in every dimension $2n+1\geq 3$ (Theorem~\ref{thm:lcs}). Thus l.c. Sasakian manifolds are exactly the conformally Sasakian manifolds of the title; we keep the name ``locally conformal'' to stress the analogy with \cite{olszak,madu}. The reason for this behaviour is that the Sasakian condition $d\eta=\Phi$, unlike the K\"ahler and almost cosymplectic conditions, is not invariant under homotheties. Unless the conformal factor is constant, these structures are neither normal nor contact metric (Corollary~\ref{cor:normal}), and they belong to the class
$
\mathcal{C}_4\oplus\mathcal{C}_5\oplus\mathcal{C}_6\oplus\mathcal{C}_{12}
$
(Remark~\ref{rem:class}). We then compute the curvature along $\xi$ and derive an integral formula showing that $\mathrm{Ric}(\xi,\xi)$ has positive mean value on compact manifolds (Proposition~\ref{prop:integral}).

Our main results concern the $k$-nullity condition, to which the $(k,\mu)'$-condition of \cite{madu} reduces (Remark~\ref{rem:convention}). In contrast with the contact metric case, a nullity function along $\xi$ is automatically constant, even in dimension $3$ (Theorem~\ref{thm:Nk}). Moreover, the condition holds if and only if $f=\exp(\sigma)$ satisfies
$
\mathrm{Hess}_g f=-(kf+C)\,g
$
for a constant $C$ (Theorem~\ref{thm:concircular}), that is, if and only if $f$ is a concircular function in the sense of Yano \cite{yano-concircular} and Tashiro \cite{tashiro}. By Yano's criterion, this means that $g$ and the underlying Sasakian metric
$
g_0=\exp(-2\sigma)g
$
have the same geodesic circles (Theorem~\ref{thm:circles}). The situation is reminiscent of Brinkmann's theorem on conformal maps between Einstein spaces \cite{brinkmann,kuhnel}. Several consequences follow (Corollaries~\ref{cor:range} and~\ref{cor:structure}):
\begin{enumerate}
	\item[-] the Lee vector field is torse-forming;
	\item[-] the level sets of the conformal factor are totally umbilical with constant mean curvature, and the metric is locally a warped product;
	\item[-] the conformal factor $\alpha=\exp(-\sigma)$ lies in an explicit interval;
	\item[-] on a compact manifold,
	$
	k=(\min\alpha)(\max\alpha)>0.
	$
\end{enumerate}
If $g$ is complete and the Lee form is not zero, then $(M,g)$ is a sphere, Euclidean space or hyperbolic space of curvature $k$ (Theorem~\ref{thm:complete}), and all three cases occur (Remark~\ref{rem:spaceforms}). Without completeness, there are examples that are not of constant curvature in every dimension $2n+1\geq 5$. They are obtained from solutions of Obata's equation on Sasakian manifolds (Proposition~\ref{prop:construction}, Example~\ref{ex:join}). The $k$-nullity condition is therefore far more rigid here than in contact metric geometry, and it differs sharply from the almost cosymplectic case, where it leads to local products \cite{madu}.

Our approach differs from that of \cite{falcitelli1}, where the $k$-nullity condition is combined with an assumption on the $\phi$-sectional curvature. Here no further curvature condition is imposed. Theorem~\ref{thm:lcs} recovers the global conformality of \cite{falcitelli1} in the Sasakian case and extends it to dimension $3$, and Theorem~\ref{thm:Nk} holds in every dimension. To the best of our knowledge, the concircular characterization (Theorems~\ref{thm:concircular} and~\ref{thm:circles}), the range of the conformal factor (Corollary~\ref{cor:range}) and the classification in the complete case (Theorem~\ref{thm:complete}) are new. In a companion paper \cite{partII}, we treat the same condition through the Sasakian metric $g_0$, where it becomes Obata's equation, and we describe the local models by means of K\"ahler cones.

The paper is organized as follows. Section~\ref{sec:sasakian} contains preliminaries, including the notions of concircular function and concircular transformation. Section~\ref{sec:lcs} characterizes l.c.\ Sasakian manifolds, and Section~\ref{sec:basic} gives the basic identities. Section~\ref{sec:nullity} studies the $k$-nullity condition, and Section~\ref{sec:consequences} its geometric consequences, together with a construction of examples, including examples that are not of constant curvature. Section~\ref{sec:conclusion} concludes with open questions.

\section{Preliminaries}\label{sec:sasakian}

Let $M$ be a $(2n + 1)$-dimensional manifold, $n\ge1$, endowed with an almost contact structure $(\phi, \xi, \eta)$, i.e., $\phi$ is a tensor field of type $(1, 1)$, $\xi$ is a vector field, and $\eta$ is a $1$-form satisfying \cite{bl}
\begin{equation}\label{equa1}
	\phi^{2} = -\mathbb{I} + \eta\otimes\xi,\quad\eta(\xi)= 1 ,\quad\eta\circ\phi = 0 \quad\mbox{and}\quad\phi\xi = 0.
\end{equation}
Then $(\phi, \xi, \eta, g)$ is called an \emph{almost contact metric structure} on $M$ if $(\phi,\xi,\eta)$ is an almost contact structure on $M$ and $g$ is a Riemannian metric on $M$ such that, for any vector fields $X$ and $Y$ on $M$,
\begin{equation}\label{equation2}
	g(\phi X, \phi Y) = g(X, Y) -  \eta(X)\,\eta(Y).
\end{equation}
It follows from \eqref{equa1} and \eqref{equation2} that $\eta=g(\xi,\cdot)$. The fundamental $2$-form is $\Phi(X,Y)=g(X,\phi Y)$. 

All manifolds are smooth and connected, and $\Gamma(\Xi)$ denotes the module of smooth sections of a vector bundle $\Xi$.

We use the conventions of \cite{bl}:
\begin{align*}
	(\vartheta\wedge\varrho)(X,Y)&=\tfrac{1}{2}\{\vartheta(X)\varrho(Y)-\vartheta(Y)\varrho(X)\},\\
	d\eta(X,Y)&=\tfrac{1}{2}\{X\eta(Y)-Y\eta(X)-\eta([X,Y])\},
\end{align*}
for which $d$ is an antiderivation with $d^2=0$ and $d(f\vartheta)=df\wedge\vartheta+f\,d\vartheta$.

The structure is \emph{contact metric} if $d\eta=\Phi$ and \emph{normal} if $[\phi,\phi]+2\,d\eta\otimes\xi=0$, where
$$
[\phi,\phi](X,Y)=\phi^{2}[X,Y]+[\phi X,\phi Y]-\phi[\phi X,Y]-\phi[X,\phi Y]
$$
is the Nijenhuis torsion of $\phi$. It is \emph{Sasakian} if
\begin{equation}\label{eq:sasakian}
	(\nabla_{X}\phi)Y =  g(X, Y)\xi - \eta(Y)X ,
\end{equation}
where $\nabla$ is the Levi-Civita connection of $g$. A Sasakian structure satisfies (see \cite{bl})
\begin{equation}\label{eq:sasakian2}
	\nabla_{X}\xi=-\phi X,\quad \mathcal{L}_{\xi}\phi=0,\quad R(X,Y)\xi=\eta(Y)X-\eta(X)Y,
\end{equation}
where $\mathcal{L}_{\xi}$ is the Lie derivative with respect to $\xi$ and $R$ is the Riemann curvature tensor of $g$. We shall use the following characterization given in \cite[Theorem~6.3]{bl}.
% NOTE (authors): check the theorem number in the edition of [bl] listed in the bibliography (2002).

\begin{theorem}\label{thm:blair63}
	An almost contact metric structure is Sasakian if and only if it is contact metric and normal.
\end{theorem}

Our curvature conventions are
$$
R(X,Y)Z=\nabla_{X}\nabla_{Y}Z-\nabla_{Y}\nabla_{X}Z-\nabla_{[X,Y]}Z,\;\; R(X,Y,Z,W)=g(R(X,Y)Z,W).
$$
For a metric $h$ we put
\begin{equation}\label{eq:G}
	G_h(X,Y,Z,W)=h(Y,Z)h(X,W)-h(X,Z)h(Y,W),
\end{equation}
for any vector fields $X$, $Y$, $Z$ and $W$ on $M$, so that $h$ has constant sectional curvature $c$ if and only if $R_h=c\,G_h$.

The Hessian of a function $f$ with respect to a metric $h$ is $\mathrm{Hess}_{h} f=\nabla^h df$, that is, $(\mathrm{Hess}_h f)(X,Y)=h(\nabla^h_{X}\mathrm{grad}_h f,Y)$; it is symmetric since $\nabla^h$ is torsion-free. A function $f$ on a Riemannian manifold $(M,h)$ is \emph{concircular} if $\mathrm{Hess}_hf=\rho\,h$ for some function $\rho$ \cite{yano-concircular,tashiro}. Its gradient is then a \emph{concircular vector field}. A vector field $V$ is \emph{torse-forming} if $\nabla_XV=\rho X+\psi(X)V$ for a function $\rho$ and a $1$-form $\psi$ \cite{yano-torse}.

The notions above go back to Yano \cite{yano-concircular} and Tashiro \cite{tashiro}, who wrote them in classical index notation. If $f_i=\partial_if$ and $f_{;ij}$ denote the components of $df$ and of its covariant derivative with respect to $h$, then
$$
f_{;ij}=\nabla_i\nabla_jf=\partial_i\partial_jf-\Gamma^k_{ij}\,\partial_kf=(\mathrm{Hess}_hf)_{ij},
$$
so the second covariant derivative of a function is exactly its Hessian.
\begin{enumerate}
	\item[(i)] Yano called a conformal change $h^*=\exp(2p)h$ a \emph{concircular transformation} if it maps geodesic circles, i.e., curves of constant first curvature and vanishing second curvature, to geodesic circles. He showed that this happens if and only if
	$$
	p_{;ij}-p_i\,p_j=\lambda\,h_{ij}
	$$
	for some function $\lambda$. In terms of the Hessian, this reads $\mathrm{Hess}_hp-dp\otimes dp=\lambda\,h$. Since
	$$
	\mathrm{Hess}_h\big(\exp(-p)\big)=-\exp(-p)\big(\mathrm{Hess}_hp-dp\otimes dp\big),
	$$
	the condition is equivalent to $\mathrm{Hess}_h(\exp(-p))=-\exp(-p)\lambda\,h$. Thus $h^*=s^{-2}h$, with $s>0$, is a concircular transformation if and only if $s$ is a concircular function of $(M,h)$. This criterion is used in Theorem~\ref{thm:circles}.
	\item[(ii)] Tashiro called a non-constant function $u$ a \emph{concircular scalar field} if
	$$
	\nabla_j\nabla_iu=\lambda\,h_{ji},\quad\text{i.e.,}\quad \mathrm{Hess}_hu=\lambda\,h,
	$$
	and its gradient $v$, which satisfies $\nabla_Xv=\lambda X$, a \emph{concircular vector field}. When $\lambda=au+b$ with constants $a$ and $b$, the field is called \emph{special concircular}. Equation~\eqref{eq:concircular} in Section~\ref{sec:consequences} is of this type, with $a=-k$ and $b=-C$.
\end{enumerate}

\section{Locally conformal Sasakian manifolds}\label{sec:lcs}

In this section we introduce the class of manifolds studied in this paper and show that it consists exactly of conformal changes of Sasakian manifolds. We follow Olszak's construction of locally conformal almost cosymplectic manifolds \cite{olszak}, with the almost cosymplectic local model replaced by a Sasakian one. If $(\phi,\xi,\eta,g)$ is an almost contact metric structure and $\sigma_t$ is a function on an open set $U_t$, then $(\phi,\xi_t,\eta_t,\exp(-2\sigma_t)g)$ is again an almost contact metric structure only if $\xi_t=\pm \exp(\sigma_t)\xi$ and $\eta_t=\pm \exp(-\sigma_t)\eta$. Indeed, $\phi^2=-\mathbb I+\eta_t\otimes\xi_t$ forces $\eta_t\otimes\xi_t=\eta\otimes\xi$, and the compatibility of the metric then fixes the factor up to sign. This is the normalization \eqref{lcsasakianmanifold} below.

\begin{definition}\label{def:lcs}
	An almost contact metric manifold $(M,\phi,\xi,\eta,g)$ is \emph{locally conformal (l.c.) Sasakian} if there exist an open cover $\{U_{t}\}_{t\in I}$ of $M$ and smooth functions $\sigma_{t}: U_{t}\to\mathbb{R}$ such that
	\begin{equation}\label{lcsasakianmanifold}
		\phi_{t} =\phi,\quad g_{t}=\exp(-2\sigma_{t})g,\quad \xi_{t} = \exp(\sigma_{t})\xi,\quad \eta_{t}=\exp(-\sigma_{t})\eta
	\end{equation}
	is a Sasakian structure on $U_{t}$.
\end{definition}

Our first result shows that the local conformal factors are determined by the structure.

\begin{theorem}\label{thm:lcs}
	An almost contact metric manifold $(M,\phi,\xi,\eta,g)$ is locally conformal Sasakian if and only if there exists a function $\sigma\in C^\infty(M)$ such that
	\begin{gather}
		d\eta = d\sigma\wedge\eta + \exp(-\sigma)\Phi, \label{eq:i}\\
		[\phi,\phi] +2\exp(-\sigma)\Phi\otimes\xi=0.	\label{eq:ii}
	\end{gather}
	In this case, the local conformal factors $\sigma_t$ of Definition~\ref{def:lcs} are the restrictions of $\sigma$ to the sets $U_t$. In particular, they agree on overlaps. The function $\sigma$ is determined by the almost contact metric structure: for every orthonormal basis $\{E_1,\ldots,E_{2n}\}$ of $\ker\eta$ at a point, and for every orthonormal $\phi$-basis $\{e_1,\ldots,e_n,\phi e_1,\ldots,\phi e_n\}$ of $\ker\eta$ at that point,
	\begin{equation}\label{eq:alpha}
		\exp(-\sigma) =-\frac{1}{2n}\sum_{a=1}^{2n}d\eta(E_{a},\phi E_{a})=-\frac{1}{n}\sum_{i=1}^{n}d\eta(e_{i},\phi e_{i}).
	\end{equation}
	Consequently, the Lee form $\omega:=d\sigma$ is exact, $(\phi,\, \exp(\sigma)\xi,\, \exp(-\sigma)\eta, \,\exp(-2\sigma)g)$ is a Sasakian structure on $M$, and $d\Phi=2\,\omega\wedge\Phi$.
\end{theorem}
\begin{proof}
	Assume that $(M,\phi,\xi,\eta,g)$ is l.c.\ Sasakian, with an open cover $\{U_t\}_{t\in I}$ and functions $\sigma_t$ as in Definition~\ref{def:lcs}. On each $U_t$ the structure \eqref{lcsasakianmanifold} is Sasakian, hence contact metric and normal by Theorem~\ref{thm:blair63}, and its fundamental $2$-form is $\Phi_t=\exp(-2\sigma_t)\Phi$. Therefore,
	$\exp(-2 \sigma_{t})\Phi   = \Phi_{t}= d\eta_{t}= d(\exp(-\sigma_{t})\eta) = \exp(-\sigma_{t}) \{ - d\sigma_{t}\wedge \eta +  d\eta \}$ and $[\phi,\phi]   = -2\, d\eta_{t}\otimes \xi_{t}= -2\exp(-\sigma_{t})\Phi\otimes\xi$, that is, on each $U_t$,
	\begin{equation}\label{eq:local}
		d\eta  = d\sigma_{t} \wedge\eta + \exp(-\sigma_{t})\Phi, \;\;
		[\phi,\phi] + 2\exp(-\sigma_{t})\Phi\otimes\xi =0.
	\end{equation}
	If $U_t\cap U_s\neq\emptyset$, subtracting the second identities in \eqref{eq:local} for $t$ and $s$ gives
	$
	\{\exp(-\sigma_{t}) - \exp(-\sigma_{s})\}\,\Phi\otimes\xi =0\quad\text{on } U_t\cap U_s .
	$
	Since $\Phi\otimes\xi$ vanishes nowhere ($\Phi(X,\phi X)=-1$ for every unit vector $X\in\ker\eta$), we get $\exp(-\sigma_t)=\exp(-\sigma_s)$ on $U_t\cap U_s$, and the injectivity of the exponential function gives $\sigma_{t}=\sigma_{s}$ there. Hence the local functions $\sigma_t$ define a unique global function $\sigma\in C^\infty(M)$ with $\sigma|_{U_t}=\sigma_t$, and \eqref{eq:local} gives \eqref{eq:i} and \eqref{eq:ii}. Conversely, let $\sigma\in C^\infty(M)$ satisfy \eqref{eq:i} and \eqref{eq:ii}, and put $\omega=d\sigma$, $\eta_{0}= \exp(-\sigma)\eta$, $\xi_{0} = \exp(\sigma)\xi$ and $g_{0} = \exp(-2\sigma)g$. Since $\eta_{0}\otimes\xi_{0}=\eta\otimes\xi$, we get $\phi^{2}= -\mathbb{I} + \eta_{0}\otimes\xi_{0}$. Also $\eta_{0}(\xi_{0})=1$ and
	$
	g_{0}(\phi X,\phi Y) =\exp(-2\sigma) \{g(X,Y) -\eta(X)\eta(Y)\}= g_{0}(X,Y)-\eta_{0}(X)\eta_{0}(Y).
	$
	Thus $(\phi,\xi_{0},\eta_{0},g_{0})$ is an almost contact metric structure on $M$, whose fundamental $2$-form is $\Phi_{0} = \exp(-2\sigma)\Phi$. By (\ref{eq:i}),
	$
	d\eta_{0}=\exp(-\sigma)\bigl(d\eta-\omega\wedge\eta\bigr)=\exp(-2\sigma) \Phi=\Phi_{0},
	$
	so this structure is contact metric. Furthermore, by \eqref{eq:ii},
	$
	[\phi,\phi]+2\,d\eta_{0}\otimes\xi_{0} = [\phi,\phi]+2\exp(-\sigma)\Phi\otimes\xi=0,
	$
	so it is normal. By Theorem~\ref{thm:blair63} it is Sasakian, and $M$ is l.c.\ Sasakian, with the cover $\{M\}$ and $\sigma_M=\sigma$. This also proves that $(\phi,\, \exp(\sigma)\xi,\, \exp(-\sigma)\eta, \,\exp(-2\sigma)g)$ is Sasakian.
	
	Next, we apply $d$ to \eqref{eq:i}. Using $d\omega=0$, $\omega\wedge\omega=0$ and $d(\exp(-\sigma))=-\exp(-\sigma)\,\omega$, we obtain
	\begin{align*}
		0&=-\omega\wedge d\eta-\exp(-\sigma)\omega\wedge\Phi+\exp(-\sigma)d\Phi\\
		&=-\omega\wedge\big(\omega\wedge\eta+\exp(-\sigma)\Phi\big)-\exp(-\sigma)\omega\wedge\Phi+\exp(-\sigma)d\Phi\\
		&=\exp(-\sigma)\{d\Phi-2\,\omega\wedge\Phi\}.
	\end{align*}
	Hence $d\Phi=2\,\omega\wedge\Phi$. Finally, let $X\in\ker\eta$ be a unit vector. Since $\eta(X)= \eta(\phi X)=0$, evaluating (\ref{eq:i}) at $(X,\phi X)$ gives
	$
	d\eta(X,\phi X) = \exp(-\sigma)g(X,\phi^2X)=-\exp(-\sigma).
	$
	Taking $X=E_a$ and summing over $a=1,\ldots,2n$ gives the first formula in \eqref{eq:alpha}. If $\{e_i,\phi e_i\}_{1\le i\le n}$ is an orthonormal $\phi$-basis, then
	$
	d\eta(\phi e_i,\phi^2e_i) =d\eta(e_i,\phi e_i),
	$
	so the sum over this basis of $2n$ vectors is twice the sum over $i=1,\ldots,n$. This gives the second formula. The right-hand side of \eqref{eq:alpha} depends only on the almost contact metric structure, which gives a second proof that the $\sigma_t$ agree on overlaps.
\end{proof}

\begin{remark}\label{rem:global}
	\emph{Theorem~\ref{thm:lcs} shows that, with the normalization \eqref{lcsasakianmanifold}, an l.c.\ Sasakian manifold is \emph{globally} conformal to a Sasakian manifold, and its Lee form $\omega=d\sigma$ is exact. Conversely, every conformal change of a Sasakian structure of the form \eqref{lcsasakianmanifold} is l.c.\ Sasakian (Example~\ref{ex:lcs}). Thus l.c.\ Sasakian manifolds are exactly the conformally Sasakian manifolds; we keep the name ``locally conformal'' to stress the analogy with \cite{olszak,madu}. This contrasts with the l.c.\ almost cosymplectic case, where the defining conditions $d\eta=0$, $d\Phi=0$ are invariant under homotheties $(\eta,g)\longmapsto(c\eta,c^2g)$, so that local conformal factors may differ by constants on overlaps \cite{olszak}. It also contrasts with locally conformal K\"ahler geometry, where the Hopf manifolds are locally but not globally conformal K\"ahler \cite{dragomir-ornea}. The reason is that the Sasakian condition $d\eta=\Phi$ is not homothety invariant. Allowing $a_{t}$-Sasakian local models ($d\eta_{t}=a_{t}\Phi_{t}$ with constants $a_{t}>0$) gives nothing new. Indeed, $(\phi,a_{t}^{-1}\xi_{t},a_{t}\eta_{t},a_{t}^{2}g_{t})$ is then Sasakian, and the homothety is absorbed by replacing $\sigma_{t}$ with $\sigma_{t}-\ln a_{t}$.}
\end{remark}

\begin{lemma}\label{lem:omegaeta}
	Let $M$ be an l.c.\ Sasakian manifold with Lee form $\omega=d\sigma$, and let $U\subset M$ be open. If $\omega\wedge\eta=0$ on $U$, then $\omega=0$ on $U$.
\end{lemma}
\begin{proof}
	On $U$ we have $\omega=\beta\eta$ with $\beta=\omega(\xi)$. By \eqref{eq:i},
	$
	0=d\omega=d\beta\wedge\eta+\beta\,d\eta=d\beta\wedge\eta+\beta \exp(-\sigma)\Phi.
	$
	Evaluating at $(X,\phi X)$ with $X\in\ker\eta$ a unit vector gives $\beta \exp(-\sigma)=0$, so $\beta=0$ and $\omega=0$ on $U$.
\end{proof}

\begin{corollary}\label{cor:normal}
	Let $(M,\phi,\xi,\eta,g)$ be an l.c.\ Sasakian manifold, and let $\sigma$ and $\omega=d\sigma$ be as in Theorem~\ref{thm:lcs}. Then:
	\begin{enumerate}
		\item[(i)] the structure is normal if and only if $\omega=0$;
		\item[(ii)] the structure is contact metric if and only if $\sigma=0$;
		\item[(iii)] in particular, if $\omega\neq 0$, the structure is neither normal nor contact metric.
	\end{enumerate}
\end{corollary}
\begin{proof}
	By \eqref{eq:i} and \eqref{eq:ii},
	$
	[\phi,\phi]+2\,d\eta\otimes\xi 	=2(\omega\wedge\eta)\otimes\xi.
	$
	If the structure is normal, then $\omega\wedge\eta=0$ on $M$, and Lemma~\ref{lem:omegaeta} gives $\omega=0$; the converse is clear. If the structure is contact metric, then $d\eta=\Phi$, and evaluating \eqref{eq:i} at $(X,\phi X)$, with $X\in\ker\eta$ a unit vector, gives $1=\exp(-\sigma)$, hence $\sigma=0$. Conversely, $\sigma=0$ gives $d\eta=\Phi$. Finally, $\omega\neq 0$ implies that $\sigma$ is not constant, so (iii) follows from (i) and (ii).
\end{proof}

The following example shows that a Sasakian manifold and an arbitrary function always produce an l.c.\ Sasakian manifold. By Theorem~\ref{thm:lcs}, every l.c.\ Sasakian manifold arises in this way.

\begin{example}\label{ex:lcs}
	\emph{Let $(N,\phi_{0},\xi_{0},\eta_{0},g_{0})$ be any $(2n+1)$-Sasakian manifold and let $\sigma\in C^\infty(N)$ be an arbitrary smooth function. On $M=N$ define
		\begin{equation}\label{eq:exstructure}
			\phi=\phi_{0},\quad \xi=\exp(-\sigma)\xi_{0},\quad \eta=\exp(\sigma)\eta_{0},\quad g=\exp(2\sigma)g_{0} .
		\end{equation}
		Then $(M,\phi,\xi,\eta,g)$ is a locally conformal Sasakian manifold. Indeed, take the cover consisting of the single open set $U=M$ and $\sigma_U=\sigma$. Then \eqref{lcsasakianmanifold} gives back $\phi_U=\phi_{0}$, $g_U=\exp(-2\sigma)g=g_{0}$, $\xi_U=\exp(\sigma)\xi=\xi_{0}$ and $\eta_U=\exp(-\sigma)\eta=\eta_{0}$, which is Sasakian. Let us check Theorem~\ref{thm:lcs} directly. Here $\omega=d\sigma$ and $\Phi= \exp(2\sigma)d\eta_{0}$. Then
		$
		d\eta=d(\exp(\sigma)\eta_{0})=\exp(\sigma)d\sigma\wedge\eta_{0} + \exp(\sigma)d\eta_{0}=\omega\wedge\eta+\exp(-\sigma)\Phi,
		$
		and, since $[\phi_{0},\phi_{0}]=-2\,d\eta_{0}\otimes\xi_{0}$,
		$
		[\phi,\phi] =-2\,d\eta_{0}\otimes\xi_{0}=-2\exp(-2\sigma)\Phi\otimes\exp(\sigma)\xi=-2\exp(-\sigma)\,\Phi\otimes\xi .
		$
		So \eqref{eq:i} and \eqref{eq:ii} hold. By Corollary~\ref{cor:normal}, the structure \eqref{eq:exstructure} is Sasakian only if $\sigma=0$.}
\end{example}

\begin{example}\label{ex:R2n+1}
	\emph{Consider the $(2n+1)$-dimensional $M =\{p\in \mathbb{R}^{2n+1}: z>0\}$, where $p=(x_1,\dots,x_n,y_1,\dots,y_n,z)$ is the standard coordinates on $\mathbb{R}^{2n+1}$. Recall the standard Sasakian structure \cite{bl}
		\begin{align*}
			\eta_{0}&=\tfrac{1}{2}\Big\{dz-\sum_{i=1}^n y_{i}\,dx_{i}\Big\},\quad \xi_{0}=2\frac{\partial}{\partial z},\quad
			g_{0}=\eta_{0}\otimes\eta_{0}+\tfrac14\sum_{i=1}^n\{dx_{i}^{2}+dy_{i}^{2}\},\\
			\phi_{0}&=\sum_{i=1}^n\Big[\Big\{\frac{\partial}{\partial x_{i}}+y_{i}\frac{\partial}{\partial z}\Big\}\otimes dy_{i}
			-\frac{\partial}{\partial y_{i}}\otimes dx_{i}\Big],
		\end{align*}
		for which $\Phi_{0}=d\eta_{0}=\frac{1}{2}\sum_{i} dx_{i}\wedge dy_{i}$. Take $\sigma=z$. By Example~\ref{ex:lcs},
		$$
		\eta=\tfrac{1}{2}\exp(z)\Big\{dz-\sum_{i} y_{i}\,dx_{i}\Big\},\quad
		\xi=2\exp(-z)\frac{\partial}{\partial z},\quad \phi=\phi_{0},
		$$
		$$
		g=\exp(2z)\Big[\tfrac14\Big\{dz-\sum_{i} y_{i}\,dx_{i}\Big\}^2+\tfrac14\sum_{i}\{dx_{i}^{2}+dy_{i}^{2}\}\Big]
		$$
		is a locally conformal Sasakian structure on $M$, with
		$$
		\omega=dz,\;\; \Phi=\tfrac{1}{2}\exp(2z)\sum_{i} dx_{i}\wedge dy_{i}.
		$$
		Directly,
		\begin{align*}
			d\eta &=\tfrac{1}{2}\exp(z)\Big\{-\sum_{i} y_{i}\,dz\wedge dx_{i}+\sum_{i} dx_{i}\wedge dy_{i}\Big\}=\omega\wedge\eta + \exp(- z)\Phi,\\
			d\Phi &= \exp(2z)\,dz\wedge\sum_{i} dx_{i}\wedge dy_{i}= 2\,\omega\wedge\Phi .
		\end{align*}
		By Corollary~\ref{cor:normal}, this structure is neither normal nor contact metric. We return to this example after Corollary~\ref{cor:range}.
	}
\end{example}

\section{Some basic identities}\label{sec:basic}

From now on, $(M,\phi,\xi,\eta,g)$ is an l.c.\ Sasakian manifold, $\sigma$ is the function of Theorem~\ref{thm:lcs}, and
$$
g_{0}=\exp(-2\sigma) g,\quad \xi_{0}= \exp(\sigma)\xi,\quad \eta_{0}=\exp(- \sigma)\eta,\quad \phi_{0}=\phi
$$
is the underlying Sasakian structure. We write
\begin{align*}
	\omega & = d\sigma,\quad g(B,X)=\omega(X),\quad \alpha= \exp(-\sigma)>0,\quad\beta=\omega(\xi),\\
	\theta & =g(\phi B,\cdot\,)=-\omega\circ\phi,\quad \nu^{2}=\|B\|^{2} - \beta^{2}.
\end{align*}
Here and below $\|\cdot\|$ is the norm with respect to $g$; norms with respect to $g_{0}$ carry a subscript. The vector field $B$ is the \emph{Lee vector field} and $\theta$ the \emph{anti-Lee form}. Note that $d\alpha=-\alpha\,\omega$, $\|\phi B\|^{2}=\nu^{2}$, $\theta(\xi)=\theta(B)=0$, and that $\theta=0$ at a point if and only if $\omega\wedge\eta=0$ there, i.e., $B=\beta\xi$. Since $\nabla\omega=\mathrm{Hess}_g\sigma$ is symmetric, $g(\nabla_{X}B,Y)=g(\nabla_{Y}B,X)$.

By Theorem~\ref{thm:lcs},
\begin{equation}\label{eq:structure}
	d\eta=\omega\wedge\eta+\alpha\,\Phi,\;\; d\Phi=2\,\omega\wedge\Phi .
\end{equation}
As in \cite{olszak}, the Levi-Civita connections $\nabla$ of $g$ and $\nabla^{0}$ of $g_{0}$ are related by
\begin{equation}\label{eq:conn}
	\nabla^{0}_{X}Y=\nabla_{X}Y-\omega(X)Y-\omega(Y)X+g(X,Y)B .
\end{equation}

% NOTE (authors): check the correspondence with the formulas of [abedi] (normalization e^{f}=e^{-2\sigma}).
The following two propositions correspond, after a change of normalization, to formulas of Abedi and Ilmakchi \cite{abedi}. We include the short proofs for completeness.

\begin{proposition}\label{prop:nablaphi}
	Let $(M, \phi,\xi,\eta,g)$ be an l.c.\ Sasakian manifold. Then
	\begin{align}
		(\nabla_{X}\phi)Y ={}&\alpha\left\{g(X,Y)\xi-\eta(Y)X\right\}+\omega(\phi Y)X-\omega(Y)\phi X\notag\\
		&-g(X,\phi Y)B + g(X,Y)\phi B,\label{eq:nablaphi}
	\end{align}
	for all $X,Y\in\Gamma(TM)$.
\end{proposition}
\begin{proof}
	By \eqref{eq:conn},
	$$
	(\nabla^{0}_{X}\phi)Y =(\nabla_{X}\phi)Y-\omega(\phi Y)X+\omega(Y)\phi X+g(X,\phi Y)B-g(X,Y)\phi B,
	$$
	and by \eqref{eq:sasakian} for the Sasakian structure $(\phi_{0},\xi_{0},\eta_{0},g_{0})$,
	$$
	(\nabla^{0}_{X}\phi)Y=g_{0}(X,Y)\xi_{0}-\eta_{0}(Y)X=\alpha\{g(X,Y)\xi-\eta(Y)X\},
	$$
	which completes the proof.
\end{proof}

\begin{proposition}\label{prop:nablaxi}
	Let $(M, \phi,\xi,\eta,g)$ be an l.c.\ Sasakian manifold. Then, for all $X\in\Gamma(TM)$,
	\begin{equation}\label{eq:nablaxi}
		\nabla_{X}\xi=-\alpha\,\phi X+\beta X-\eta(X)B .
	\end{equation}
	In particular $\nabla_\xi\xi=-B+\beta\xi$ and $\nabla_B\xi=-\alpha\,\phi B$. Moreover, with $\tilde hX:=\nabla_{X}\xi-\beta X+\eta(X)B$ and $h:=\frac{1}{2}\mathcal{L}_\xi\phi$,
	\begin{align}\label{eq:h}
		\tilde h &=-\alpha\,\phi,\quad h=-\tfrac{1}{2}\,\theta\otimes\xi,\quad\phi\circ h   =0,\nonumber\\
		h\circ\phi&=-\tfrac{1}{2}(\omega-\beta\eta)\otimes\xi, \quad \mathcal{L}_\xi\eta=\beta\,\eta-\omega.
	\end{align}
\end{proposition}
\begin{proof}
	By \eqref{eq:conn},
	$$
	\nabla^{0}_{X}\xi_{0} = \exp(\sigma)\{\omega(X)\xi+\nabla^{0}_{X}\xi\}=\exp(\sigma)\{\nabla_{X}\xi-\beta X+\eta(X)B\} =\exp(\sigma)\tilde hX.
	$$
	Since $\nabla^{0}\xi_{0}=-\phi$ by \eqref{eq:sasakian2}, $\tilde h=-\alpha\phi$, which is \eqref{eq:nablaxi}. For any function $f$,
	$$
	(\mathcal{L}_{f\xi}\phi)X=[f\xi,\phi X]-\phi[f\xi,X]=f(\mathcal{L}_\xi\phi)X-((\phi X)f)\,\xi .
	$$
	With $f=\exp(\sigma)$ and $\mathcal L_{\xi_{0}}\phi=0$ this gives $0=\exp(\sigma)\{2hX-\omega(\phi X)\xi\}$, i.e., $hX=-\frac{1}{2}\theta(X)\xi$. Then $\phi\circ h=0$, while by \eqref{equation2}
	$$
	h(\phi X)=-\tfrac12\theta(\phi X)\xi=-\tfrac12 g(\phi B,\phi X)\xi=-\tfrac12\{\omega(X)-\beta\eta(X)\}\xi.
	$$
	Finally, $(\mathcal{L}_\xi\eta)Y=2\,d\eta(\xi,Y)$, since $\eta(\xi)=1$, and \eqref{eq:structure} gives $2\,d\eta(\xi,Y)=\beta\eta(Y)-\omega(Y)$, and the proof is complete.
\end{proof}
\begin{remark}\label{rem:class}
	\emph{From \eqref{eq:nablaphi}, writing $\psi=\omega\circ\phi=-\theta$, one obtains
		\begin{align*}
			(\nabla_{X}\Phi)(Y,Z)={}&\alpha\{g(X,Z)\eta(Y)-g(X,Y)\eta(Z)\}\\
			&-\beta\{\Phi(X,Z)\eta(Y)-\Phi(X,Y)\eta(Z)\}\\
			& +g(\phi X,\phi Y)\psi(Z)-g(\phi X,\phi Z)\psi(Y)\\
			&-\Phi(X,Y)\psi(\phi Z)+\Phi(X,Z)\psi(\phi Y)\\
			&+\eta(X)\{\eta(Z)\theta(Y)-\eta(Y)\theta(Z)\},
		\end{align*}
		where $(\nabla_{X}\Phi)(Y,Z)=g(Y,(\nabla_{X}\phi)Z)$ and $(\nabla_\xi\eta)(\phi Y)=\theta(Y)$. The four groups of terms have the form of the defining relations of the classes $\mathcal{C}_6$, $\mathcal{C}_5$, $\mathcal{C}_4$ and $\mathcal{C}_{12}$ of Chinea and Gonz\'alez \cite{chinea-gonzalez}; the $\mathcal{C}_4$-part vanishes identically when $n=1$. Hence l.c.\ Sasakian manifolds belong to $\mathcal{C}_4\oplus\mathcal{C}_5\oplus\mathcal{C}_6\oplus\mathcal{C}_{12}$, with nowhere vanishing $\mathcal{C}_6$-component $\alpha$. The $\mathcal{C}_5$-component is determined by $\beta$, and the $\mathcal{C}_4$- and $\mathcal{C}_{12}$-components by $\theta$,  see \cite{chinea-marrero,falcitelli1} for more details.}
\end{remark}
\begin{remark}\label{rem:reeb}
	\emph{By \eqref{eq:structure}, $\eta\wedge(d\eta)^n=\alpha^n\,\eta\wedge\Phi^n\neq0$. So $\eta$ is a contact form and $\ker\eta$ is not integrable, whereas it is integrable in the l.c.\ almost cosymplectic case. Since $2\,d\eta(\xi,\cdot)=\beta\eta-\omega$, the characteristic vector field $\xi$ is the Reeb field of $\eta$ exactly where $\theta=0$. By Lemma~\ref{lem:omegaeta}, this can hold on an open set only if $\omega=0$ there.}
\end{remark}

\begin{lemma}\label{lem:dbeta}
	Let $(M, \phi,\xi,\eta,g)$ be an l.c.\ Sasakian manifold. Then
	$$
	g(\nabla_{X}B,\xi)=X(\beta)-\alpha\,\theta(X)-\beta\,\omega(X)+\eta(X)\|B\|^{2}
	$$
	for all $X \in\Gamma(T M)$. Consequently,
	\begin{equation}\label{eq:nablaxiB}
		\nabla_\xi B=\mathrm{grad}\,\beta-\alpha\,\phi B-\beta B+\|B\|^{2}\xi,\;\;
		\eta(\nabla_\xi B)=\xi(\beta)-\beta^{2}+\|B\|^{2} .
	\end{equation}
\end{lemma}
\begin{proof}
	We have $X(\beta)=g(\nabla_{X}B,\xi)+g(B,\nabla_{X}\xi)$ and, by \eqref{eq:nablaxi},
	$$
	g(B,\nabla_{X}\xi)=\alpha\,\theta(X)+\beta\omega(X)-\eta(X)\|B\|^{2}.
	$$
	Since $\nabla\omega$ is symmetric, $g(\nabla_\xi B,X)=g(\nabla_{X}B,\xi)$.
\end{proof}

\begin{theorem}\label{thm:Rxi}
	Let $(M, \phi,\xi,\eta,g)$ be an l.c.\ Sasakian manifold. Then
	\begin{align}
		R(X,Y)\xi={}&\alpha^{2}\left\{\eta(Y)X-\eta(X)Y\right\}
		+\big\{\eta(\nabla_{X}B)+\beta\,\omega(X)-\eta(X)\|B\|^{2}\big\}Y \notag\\
		&-\big\{\eta(\nabla_{Y}B)+\beta\,\omega(Y)-\eta(Y)\|B\|^{2}\big\}X \notag\\
		&-\left\{\omega(X)\eta(Y)-\omega(Y)\eta(X)\right\}B
		+\eta(X)\nabla_{Y}B-\eta(Y)\nabla_{X}B, \label{eq:Rxi}
	\end{align}
	for all $X,Y\in\Gamma(TM)$.
\end{theorem}
\begin{proof}
	Differentiating \eqref{eq:nablaxi} and using $X(\alpha)=-\alpha\omega(X)$,
	$$
	\nabla_{X}\nabla_{Y}\xi =\alpha\omega(X)\phi Y-\alpha\nabla_{X}(\phi Y) + X(\beta)Y +\beta\nabla_{X}Y
	-X(\eta(Y))B-\eta(Y)\nabla_{X}B.
	$$
	Antisymmetrizing in $X,Y$, subtracting $\nabla_{[X,Y]}\xi=-\alpha\phi[X,Y]+\beta[X,Y]-\eta([X,Y])B$ and using $X(\eta(Y))-Y(\eta(X))-\eta([X,Y])=2\,d\eta(X,Y)$, we get
	\begin{align*}
		R(X,Y)\xi &=  \alpha\left\{\omega(X)\phi Y-\omega(Y)\phi X\right\}
		-\alpha\left\{(\nabla_{X}\phi)Y-(\nabla_{Y}\phi)X\right\}\\
		&\quad +X(\beta)Y-Y(\beta)X -2\,d\eta(X,Y)B-\eta(Y)\nabla_{X}B+\eta(X)\nabla_{Y}B.
	\end{align*}
	By \eqref{eq:nablaphi}, and since $g(X,\phi Y)-g(Y,\phi X)=2\Phi(X,Y)$,
	\begin{align*}
		(\nabla_{X}\phi)Y-(\nabla_{Y}\phi)X & =\alpha\left\{\eta(X)Y-\eta(Y)X\right\}+\omega(\phi Y)X -\omega(\phi X)Y\\
		& \quad +\omega(X)\phi Y-\omega(Y)\phi X - 2\Phi(X,Y)B.
	\end{align*}
	After substitution the terms $\alpha\{\omega(X)\phi Y-\omega(Y)\phi X\}$ cancel, and by \eqref{eq:structure} the coefficient of $B$ becomes $2(\alpha\Phi-d\eta)(X,Y)=-\{\omega(X)\eta(Y)-\omega(Y)\eta(X)\}$. Finally, insert $X(\beta)=\eta(\nabla_{X}B)+\alpha\theta(X)+\beta\omega(X)-\eta(X)\|B\|^2$ from Lemma~\ref{lem:dbeta}; its term $\alpha\theta(X)Y=-\alpha\,\omega(\phi X)Y$ cancels the remaining term $\alpha\,\omega(\phi X)Y$, and likewise for $Y$.
\end{proof}

\begin{theorem}\label{thm:jacobi}
	Let $(M, \phi,\xi,\eta,g)$ be an l.c.\ Sasakian manifold. Then the Jacobi operator $\ell=R(\cdot,\xi)\xi$ satisfies
	\begin{align}
		\ell X={}&\alpha^{2}\left\{X-\eta(X)\xi\right\}-\xi(\beta)X+\left\{X(\beta)-\alpha\theta(X)\right\}\xi
		-\omega(X)B-\nabla_{X}B \notag\\
		&+\eta(X)\left\{\mathrm{grad}\,\beta-\alpha\phi B+\|B\|^{2}\xi\right\}, \label{eq:ell}\\
		\mathrm{Ric}(\xi,\xi)={}&\mathrm{trace}\,\ell=2n\,\alpha^{2}-(2n-1)\,\xi(\beta)-\mathrm{div}\,B .
		\label{eq:traceell}
	\end{align}
\end{theorem}
\begin{proof}
	Put $Y=\xi$ in \eqref{eq:Rxi}. By Lemma~\ref{lem:dbeta}, $\eta(\nabla_{X}B)+\beta\omega(X)-\eta(X)\|B\|^{2}=X(\beta)-\alpha\theta(X)$ and $\eta(\nabla_\xi B)+\beta^{2}-\|B\|^{2}=\xi(\beta)$; substituting these and \eqref{eq:nablaxiB} gives \eqref{eq:ell}. For the trace, use $\mathrm{trace}(X\mapsto X-\eta(X)\xi)=2n$, $\theta(\xi)=0$, $\eta(\phi B)=0$, $\mathrm{trace}(\omega\otimes B)=\|B\|^{2}$ and $\mathrm{trace}\,\nabla B=\mathrm{div}\,B$.
\end{proof}

For $\omega=0$ and $\alpha=1$, \eqref{eq:traceell} reduces to the Sasakian identity $\mathrm{Ric}(\xi,\xi)=2n$. The trace formula has the following global consequence.

\begin{proposition}\label{prop:integral}
	On an l.c.\ Sasakian manifold, $\mathrm{div}\,\xi=2n\beta$ and
	\begin{equation}\label{eq:ricxi}
		\mathrm{Ric}(\xi,\xi)=2n\alpha^2+2n(2n-1)\beta^2-(2n-1)\,\mathrm{div}(\beta\xi)-\mathrm{div}\,B .
	\end{equation}
	In particular, if $M$ is compact, then
	\begin{equation}\label{eq:integral}
		\int_M\mathrm{Ric}(\xi,\xi)\,dV_g=2n\int_M\big(\alpha^2+(2n-1)\beta^2\big)\,dV_g>0 .
	\end{equation}
	Consequently, no compact l.c.\ Sasakian manifold satisfies $\mathrm{Ric}(\xi,\xi)\le0$ everywhere.
\end{proposition}
\begin{proof}
	By \eqref{eq:nablaxi}, $\mathrm{div}\,\xi=\mathrm{trace}(-\alpha\phi+\beta\,\mathbb I-\eta\otimes B)=(2n+1)\beta-\eta(B)=2n\beta$. Hence $\xi(\beta)=\mathrm{div}(\beta\xi)-\beta\,\mathrm{div}\,\xi=\mathrm{div}(\beta\xi)-2n\beta^2$, and \eqref{eq:traceell} becomes \eqref{eq:ricxi}. Integrating and using the divergence theorem gives \eqref{eq:integral}, which is positive since $\alpha>0$.
\end{proof}

\section{The $k$-nullity condition}\label{sec:nullity}

\begin{definition}
	An l.c.\ Sasakian manifold is an \emph{l.c.\ Sasakian $N(k)$-manifold} if $\xi$ belongs to the $k$-nullity distribution of $g$ in the sense of Tanno \cite{tanno88} (see also \cite{bkp}), i.e., if for a constant $k$
	\begin{equation}\label{eq:Nk}
		R(X,Y)\xi=k\{\eta(Y)X-\eta(X)Y\}.
	\end{equation}
	We refer to \eqref{eq:Nk} as the \emph{$k$-nullity condition}.
\end{definition}

Then $\ell=-k\phi^{2}$, $\mathrm{Ric}(\xi,\xi)=2nk$, and every plane containing $\xi$ has sectional curvature $k$. We write $\gamma:=\xi(\beta)=\xi(\omega(\xi))$. The following theorem shows in particular that the constancy of $k$ need not be assumed.
 \begin{theorem}\label{thm:Nk}
	Let $M$ be a $(2n+1)$-dimensional l.c.\ Sasakian manifold and $K\in C^{\infty}(M)$ such that
	\begin{equation}\label{eq:NK}
		R(X,Y)\xi= K\{\eta(Y)X-\eta(X)Y\}.
	\end{equation}
	Then $K$ is locally constant. Moreover, on every connected component, where $K=k$ is constant,
	\begin{equation}\label{eq:hessomega}
		\nabla\omega=(\gamma+\|B\|^{2})\,g-\omega\otimes\omega,
		\;\;
		k=\alpha^{2}-2\gamma-\|B\|^{2} ,
	\end{equation}
	and, for all $X$,
	\begin{align}
		\nabla_{X}B & =(\gamma+\|B\|^{2})X-\omega(X)B, \label{eq:nablaB}\\
		\nabla_\xi B & =(\gamma+\|B\|^{2})\xi-\beta B,\;\; \nabla_BB=\gamma B, \label{eq:nablaBB}\\
		d\beta&=\gamma\,\eta+\alpha\,\theta,\quad\text{i.e.,}\quad \mathrm{grad}\,\beta=\gamma\,\xi+\alpha\,\phi B, \label{eq:dbeta}\\
		\nabla_{X}(\phi B)&=\alpha\,\omega(X)\xi-\alpha\beta X+\gamma\,\phi X-\theta(X)B, \label{eq:nablaphiB}\\
		d\gamma&=-(\alpha^{2}+\gamma)\,\omega . \label{eq:dgamma}
	\end{align}
	Furthermore, there is a constant $C$ such that $\gamma=\alpha^{2}+C\alpha$ and $\|B\|^{2}=-k-\alpha^{2}-2C\alpha$. Conversely, if the first identity in \eqref{eq:hessomega} holds and $\alpha^{2}-2\gamma-\|B\|^{2}= k$ is constant, then \eqref{eq:Nk} holds.
\end{theorem}
\begin{proof}
	All terms in \eqref{eq:Rxi} are tensorial in $X$ and $Y$, so we may evaluate it at a point on vectors $X,Y\in\mathcal{D}_1:=\ker\eta$. Since $\eta(X)=\eta(Y)=0$, the terms $\alpha^{2}\{\eta(Y)X-\eta(X)Y\}$, $\{\omega(X)\eta(Y)-\omega(Y)\eta(X)\}B$ and $\eta(X)\nabla_YB-\eta(Y)\nabla_XB$ vanish, as do the terms containing $\|B\|^2$, and \eqref{eq:Rxi} reduces to
	$
	R(X,Y)\xi=\{\eta(\nabla_{X}B)+\beta\omega(X)\}Y-\{\eta(\nabla_{Y}B)+\beta\omega(Y)\}X .
	$
	The right-hand side of \eqref{eq:NK} vanishes for $X,Y\in\mathcal D_1$. Hence
	$
	\{\eta(\nabla_{X}B)+\beta\omega(X)\}Y-\{\eta(\nabla_{Y}B)+\beta\omega(Y)\}X=0 .
	$
	Since $\mathrm{rank}\,\mathcal{D}_1=2n\ge2$, every nonzero $X\in\mathcal D_1$ can be paired with some $Y\in\mathcal D_1$ linearly independent of it, and both coefficients must vanish. Thus $\eta(\nabla_{X}B)=-\beta\omega(X)$ for all $X\in\mathcal{D}_1$. By Lemma~\ref{lem:dbeta}, $X(\beta)=\alpha\theta(X)$ on $\mathcal{D}_{1}$. As $\theta(\xi)=0$, this is \eqref{eq:dbeta}, and inserting \eqref{eq:dbeta} into \eqref{eq:nablaxiB} gives the first formula in \eqref{eq:nablaBB}.
	
	Put $Y=\xi$ in \eqref{eq:Rxi} and use \eqref{eq:NK}. In view of the above and Lemma~\ref{lem:dbeta}, $\eta(\nabla_{X}B)+\beta\omega(X)-\eta(X)\|B\|^{2}=\gamma\eta(X)$ and $\eta(\nabla_\xi B)+\beta^{2}-\|B\|^{2}=\gamma$, hence
	\begin{equation}\label{eq:nablaB0}
		\nabla_{X}B=(\alpha^{2}-K)\left\{X-\eta(X)\xi\right\}-\gamma X+(2\gamma+\|B\|^{2})\eta(X)\xi-\omega(X)B .
	\end{equation}
	By \eqref{eq:nablaphi} with $Y=B$ and $\omega(\phi B)=0$,
	$
	(\nabla_{X}\phi)B=\alpha\,\omega(X)\xi-\alpha\beta X-\|B\|^2\phi X-\theta(X)B+\omega(X)\phi B,
	$
	and by \eqref{eq:nablaB0}, $\phi\nabla_{X}B=(\alpha^2-K-\gamma)\phi X-\omega(X)\phi B$. Adding these two expressions gives
	$\nabla_{X}(\phi B)=\alpha\omega(X)\xi-\alpha\beta X+c\,\phi X-\theta(X)B$, where
	$c:=\alpha^{2}-K-\gamma-\|B\|^{2}$. Since $2\,d\theta(X,Y)=g(\nabla_{X}(\phi B),Y)-g(\nabla_{Y}(\phi B),X)$ and $g(\phi X,Y)=-\Phi(X,Y)$, this gives $d\theta=\alpha\,\omega\wedge\eta-c\,\Phi+\omega\wedge\theta$. Using \eqref{eq:dbeta}, \eqref{eq:structure} and $d\alpha=-\alpha\omega$,
	$
	0=d(d\beta)=d(\gamma\eta+\alpha\theta)=\big(d\gamma+(\gamma+\alpha^2)\omega\big)\wedge\eta+\alpha(\gamma-c)\,\Phi .
	$
	Evaluating at $(X,\phi X)$, where $X\in\mathcal D_{1}$ is a unit vector, we see that the first term vanishes and $\Phi(X,\phi X)=-1$. Hence $c=\gamma$, i.e., $K=\alpha^{2}-2\gamma-\|B\|^{2}$. Then $(d\gamma+(\gamma+\alpha^2)\omega)\wedge\eta=0$, so that
	\begin{equation}\label{eq:f}
		d\gamma+(\gamma+\alpha^2)\omega=f_{0}\,\eta
	\end{equation}
	for some function $f_{0}$. Substituting $\alpha^{2}-K =2\gamma+\|B\|^{2}$ into \eqref{eq:nablaB0} and into the formula for $\nabla_{X}(\phi B)$ gives \eqref{eq:nablaB}, the first identity in \eqref{eq:hessomega}, and \eqref{eq:nablaphiB}; taking $X=B$ in \eqref{eq:nablaB} gives $\nabla_{B}B=\gamma B$. 
	
	By \eqref{eq:nablaB}, $d\|B\|^{2}=2g(\nabla B,B)=2\gamma\,\omega$. Hence, by \eqref{eq:f},
	$dK=-2\alpha^{2}\omega-2d\gamma-2\gamma\omega=-2f_{0}\eta$. Then $0=d(dK)=-2(df_{0}\wedge\eta + f_{0}\,d\eta)$, and evaluating at $(X,\phi X)$ as above gives $f_{0}\alpha=0$. Thus $f_{0}=0$, $K$ is locally constant, and \eqref{eq:f} becomes \eqref{eq:dgamma}. Since $\alpha=\exp(-\sigma)$ and $\omega=d\sigma$, \eqref{eq:dgamma} says $d(\exp(\sigma)\gamma)=d(\exp(-\sigma))$, which gives $\gamma=\alpha^{2}+C\alpha$ on each component, and the expression for $\|B\|^2$ follows from \eqref{eq:hessomega}.
	
	\emph{Converse.} Assume the first identity in \eqref{eq:hessomega}. Then $\eta(\nabla_{X}B)+\beta\omega(X)-\eta(X)\|B\|^{2}=\gamma\eta(X)$ and
	\begin{align*}
		\eta(X)\nabla_{Y}B -\eta(Y)\nabla_{X}B   ={}&(\gamma+\|B\|^{2})\left\{\eta(X)Y-\eta(Y)X\right\}\\
		&-\left\{ \eta(X)\omega(Y)-\eta(Y)\omega(X)\right\}B.
	\end{align*}
	Substituting into \eqref{eq:Rxi}, the $B$-terms cancel and
	\begin{equation}\label{eq:Kfunction}
		R(X,Y)\xi =(\alpha^{2}-2\gamma-\|B\|^{2})\{\eta(Y)X-\eta(X)Y\}.
	\end{equation}
	This completes the proof.
\end{proof}

\begin{remark}\label{rem:schur}
	\emph{Theorem~\ref{thm:Nk} is a Schur-type result. A nullity function along $\xi$ is automatically constant, in \emph{every} dimension, including $\dim M=3$. This is in contrast with contact metric geometry, where generalized $(\kappa,\mu)$-spaces with non-constant $\kappa$ and $\mu$ exist in dimension $3$ \cite{koufogiorgos-tsichlias}. The reason is the term $\alpha\,\Phi$ in $d\eta$: the identity $d(d\beta)=0$ separates into a component along $\Phi$ and a component in the ideal generated by $\eta$, and these determine $K$ and its derivative, respectively. Note also that the derivation of \eqref{eq:Kfunction} in the converse part uses only the first identity in \eqref{eq:hessomega}. This is used in Theorem~\ref{thm:circles}.}
\end{remark}
\begin{remark}\label{rem:convention}
	\emph{In \cite[(3.1)]{madu} the $(k,\mu)'$-condition
		$$
		R(X,Y)\xi=k\{\eta(Y)X-\eta(X)Y\}+\mu\{\eta(Y)h'X-\eta(X)h'Y\}
		$$
		is considered. Here it reduces to \eqref{eq:Nk}, whatever convention is used for $h'$.
		\begin{enumerate}
			\item[(i)] If $h'=\phi\circ h$, then $h'=0$ by Proposition~\ref{prop:nablaxi}.
			\item[(ii)] If $h'=h\circ\phi$, then $h'=-\frac12(\omega-\beta\eta)\otimes\xi$. Pairing the condition with $\xi$ and using $g(R(X,Y)\xi,\xi)=0$ gives $\mu\,\omega\wedge\eta=0$. On the open set where $\mu\neq0$, Lemma~\ref{lem:omegaeta} gives $\omega=0$, hence $h'=0$ there. (In the almost cosymplectic case the two conventions differ only by a sign, since $h$ anticommutes with $\phi$; here they do not.)
			\item[(iii)] More generally, if $T$ is a $(1,1)$-tensor field whose restriction to $\ker\eta$ is skew-symmetric, then the condition
			$$
			R(X,Y)\xi=k\{\eta(Y)X-\eta(X)Y\}+\mu\{\eta(Y)TX-\eta(X)TY\}
			$$
			forces $\mu T=0$ on $\ker\eta$, because the Jacobi operator $\ell$ is symmetric and satisfies $g(\ell X,\xi)=0$. This applies to the natural candidate $T=\tilde h=-\alpha\phi$.
		\end{enumerate}
		This is why we work directly with the $k$-nullity condition.}
\end{remark}
\begin{remark}\label{rem:Nk}
	\emph{For $\omega=0$ ($\sigma$ constant) one gets $k=\alpha^{2}=\exp(-2\sigma)$, which corresponds to a homothetic Sasakian structure. In general $k=\alpha^{2}-2\xi(\omega(\xi))-\|B\|^{2}$, which is the counterpart of the relation between $k$, $\mu$ and $\omega(\xi)$ in \cite[Theorem~3.4]{madu}.
	}
\end{remark}

\section{Geometric consequences of the $k$-nullity condition}\label{sec:consequences}

In this section $M$ is a connected l.c.\ Sasakian manifold with $\omega\neq 0$, i.e., $\sigma$ is not constant, and
$$
f:=\exp(\sigma)=\alpha^{-1}
$$
is the inverse conformal factor, so that $df=f\omega$, $\mathrm{grad}_gf=fB$ and $\|df\|=f\|B\|$. Since $\nabla df=df\otimes\omega+f\nabla\omega$, we have
\begin{equation}\label{eq:hessf}
	\mathrm{Hess}_gf=f\,(\nabla\omega+\omega\otimes\omega).
\end{equation}

\subsection{A concircular characterization}

\begin{theorem}\label{thm:concircular}
	Let $M$ be a connected l.c.\ Sasakian manifold with $\omega\neq 0$, and let $k\in\mathbb{R}$. Then $M$ satisfies the $k$-nullity condition \eqref{eq:Nk} if and only if there is a constant $C$ such that
	\begin{equation}\label{eq:concircular}
		\mathrm{Hess}_g f=-(kf+C)\,g .
	\end{equation}
	In this case $C$ is the constant of Theorem~\ref{thm:Nk}, and
	\begin{equation}\label{eq:firstint}
		\|df\|^2+kf^2+2Cf+1=0 .
	\end{equation}
\end{theorem}
\begin{proof}
	($\Longrightarrow$) By Theorem~\ref{thm:Nk}, $\nabla\omega+\omega\otimes\omega=(\gamma+\|B\|^2)g$ with $\gamma+\|B\|^2=(\alpha^2+C\alpha)+(-k-\alpha^2-2C\alpha)=-k-C\alpha$. By \eqref{eq:hessf} and $f\alpha=1$, $\mathrm{Hess}_gf=-f(k+C\alpha)g=-(kf+C)g$. Moreover, $\|df\|^2=f^2\|B\|^2=f^2(-k-\alpha^2-2C\alpha)=-kf^2-1-2Cf$.
	
	($\Longleftarrow$) By \eqref{eq:hessf} and \eqref{eq:concircular}, $\nabla\omega=\lambda g-\omega\otimes\omega$ with $\lambda:=-k-C\alpha$. Since $\nabla_\xi\xi=-B+\beta\xi$,
	$
	\gamma=\xi(\beta)=(\nabla_\xi\omega)(\xi)+\omega(\nabla_\xi\xi)=(\lambda-\beta^2)+(\beta^2-\|B\|^2)=\lambda-\|B\|^2 ,
	$
	so the first identity in \eqref{eq:hessomega} holds, and \eqref{eq:Kfunction} gives $R(X,Y)\xi=K\{\eta(Y)X-\eta(X)Y\}$ with $K:=\alpha^2-2\gamma-\|B\|^2=\alpha^2-2\lambda+\|B\|^2$. By Theorem~\ref{thm:Nk}, $K$ is constant. On the other hand, $Q:=\|df\|^2+kf^2+2Cf$ is constant, since by \eqref{eq:concircular}
	$
	dQ=2\,\mathrm{Hess}_gf(\mathrm{grad}_gf,\cdot)+2(kf+C)\,df=0 .
	$
	Multiplying the expression of $K$ by $f^2$ and using $f\alpha=1$ and $f^2\|B\|^2=\|df\|^2$,
	$Kf^2=1+2(kf+C)f+\|df\|^2=1+kf^2+Q$, i.e., $(K-k)f^2=1+Q$. 
	Since $\omega\neq 0$, the positive function $f$ is not constant, hence neither is $f^2$. Therefore $K=k$ and $Q=-1$, which are \eqref{eq:Nk} and \eqref{eq:firstint}. Finally, comparing $\gamma+\|B\|^2=\lambda=-k-C\alpha$ with the expression $-k-C'\alpha$ given by Theorem~\ref{thm:Nk} shows $C'=C$, since $\alpha$ is not constant.
\end{proof}

Thus the $k$-nullity condition, a condition on the curvature in one direction, is equivalent to a single linear second-order equation for the conformal factor. In fact, the factor in \eqref{eq:concircular} need not be prescribed, and Yano's criterion gives the condition a purely metric meaning.

\begin{theorem}\label{thm:circles}
	Let $M$ be a connected l.c.\ Sasakian manifold with $\omega\neq 0$, and let $g_{0}=\exp(-2\sigma)g$ be its Sasakian metric. The following conditions are equivalent:
	\begin{enumerate}
		\item[(a)] $R(X,Y)\xi=K\{\eta(Y)X-\eta(X)Y\}$ for some function $K\in C^\infty(M)$;
		\item[(b)] $M$ is an l.c.\ Sasakian $N(k)$-manifold for some constant $k$;
		\item[(c)] $f=\exp(\sigma)$ is a concircular function of $(M,g)$, i.e., $\mathrm{Hess}_gf=\rho\,g$ for some function $\rho$;
		\item[(d)] the conformal change $g\longmapsto g_{0}$ is a concircular transformation, i.e., $g$ and $g_{0}$ have the same geodesic circles.
	\end{enumerate}
	In this case $\rho=-(kf+C)$, with $k$ and $C$ as in Theorem~\ref{thm:concircular}.
\end{theorem}
\begin{proof}
	(a)$\Longrightarrow$(b) is Theorem~\ref{thm:Nk}, since $M$ is connected, and (b)$\Longrightarrow$(a) is trivial. (b)$\Longrightarrow$(c) and the value of $\rho$ follow from Theorem~\ref{thm:concircular}. (c)$\Longrightarrow$(a): if $\mathrm{Hess}_gf=\rho\,g$, then by \eqref{eq:hessf}, $\nabla\omega=\lambda g-\omega\otimes\omega$ with $\lambda:=\alpha\rho$. As in the proof of Theorem~\ref{thm:concircular}, $\gamma=\lambda-\|B\|^2$, so the first identity in \eqref{eq:hessomega} holds, and \eqref{eq:Kfunction} gives (a) with $K=\alpha^2-2\gamma-\|B\|^2$. Finally, (c)$\Longleftrightarrow$(d) is Yano's criterion recalled in Section~\ref{sec:sasakian}(i), applied with $h=g$, $s=f$ and $h^*=f^{-2}g=g_{0}$.
\end{proof}
Condition (d) involves neither the curvature nor the vector field $\xi$. The $k$-nullity condition says exactly that the Sasakian metric $g_{0}$ is concircularly related to $g$. This is reminiscent of Brinkmann's theorem on conformal mappings between Einstein spaces \cite{brinkmann,kuhnel}, where concircular functions play the same central role.

\subsection{The range of the conformal factor}

\begin{corollary}\label{cor:range}
	Let $M$ be a connected l.c.\ Sasakian $N(k)$-manifold with $\omega\neq 0$, and let $C$ be the constant in Theorem~\ref{thm:Nk}. Then $C^2>k$ and $\alpha^2+2C\alpha+k=-\|B\|^2\le 0$. Consequently, $-C-\sqrt{C^2-k}\ \le\ \alpha\ \le\ -C+\sqrt{C^2-k}$, 
	with equality precisely at points where $B=0$. If $M$ is compact, then $\min_M\alpha=-C-\sqrt{C^2-k}$, $\max_M\alpha=-C+\sqrt{C^2-k}$, 
	\begin{equation}\label{eq:kcompact}
		k=\big(\min_{M}\alpha\big)\big(\max_{M}\alpha\big)>0,
	\end{equation}
	and $C=-\frac12(\min_M\alpha+\max_M\alpha)<0$.
	Moreover, using Proposition \ref{prop:integral} and $\mathrm{Ric}(\xi,\xi)=2nk$,
	$$
	k=\frac{1}{\mathrm{Vol}(M,g)} \int_M\bigl[\alpha^2+(2n-1)\beta^2\bigr]\,dV_g.
	$$
\end{corollary}
\begin{proof}
	The identity  $\|B\|^{2}=-k-\alpha^{2}-2C\alpha$ from Theorem~\ref{thm:Nk}, since $\omega$ is not identically zero, $\|B\|^2$ is positive somewhere. Hence the quadratic $q(t)=t^2+2Ct+k$ takes negative values, so $C^2>k$ and $q(\alpha)\le 0$ gives the stated interval. If $M$ is compact, $\alpha$ attains its extrema. At an extremum $d\alpha =  - \alpha \omega=0$, hence $B=0$ and $q(\alpha)=0$.  Since $\alpha$ is not constant, $\min\alpha$ and $\max\alpha$ are the two distinct roots of $q$. Their product is $k$ and their sum is $-2C$. Both roots are positive, so $k>0$ and $C<0$. The integral identity follows from Proposition \ref{prop:integral} and the $k$-nullity condition.
\end{proof}
Corollary~\ref{cor:range} combines a pointwise identity with a global one. The product of the extreme values of the conformal factor equals the mean value of $\alpha^2+(2n-1)\beta^2$. Theorem~\ref{thm:Nk} also shows that the $k$-nullity condition is rarely satisfied, since $\|B\|^2=-k-\alpha^2-2C\alpha$ must be a function of $\alpha$. For instance, For the structure of Example~\ref{ex:R2n+1}, the vector fields $\xi_0=2\partial_z$, $E_i=2(\partial_{x_i}+y_i\partial_z)$ and $E_{n+i}=2\partial_{y_i}$ form a $g_0$-orthonormal frame, and $dz(\xi_0)=2$, $dz(E_i)=2y_i$, $dz(E_{n+i})=0$. Hence $\|dz\|^2_{g_0}=4\big(1+\sum_iy_i^2\big)$ and
$
\|B\|^2=\exp(-2z)\|dz\|^2_{g_0}=4\exp(-2z)\Big(1+\sum_{i=1}^ny_i^2\Big),
$
which is not a function of $\alpha=\exp(-z)$. Hence this structure satisfies the $k$-nullity condition for no $k$, on any open subset of $\mathbb{R}^{2n+1}$. 

\subsection{The local geometry}

\begin{corollary}\label{cor:structure}
	Let $M$ be a connected l.c.\ Sasakian $N(k)$-manifold with $\omega\neq 0$, and let $f$ and $C$ be as in Theorem~\ref{thm:concircular}.
	\begin{enumerate}
		\item[(a)] The vector field $\mathrm{grad}_gf=\exp(\sigma)B$ is concircular, $\nabla_X(\exp(\sigma)B)=-(kf+C)X$, and the Lee vector field is torse-forming,
		$\nabla_XB=-(k+C\alpha)X-\omega(X)B$.
		\item[(b)] Every plane containing $\xi$ has sectional curvature $k$.
		\item[(c)] The zeros of $B$, i.e., the critical points of $\sigma$, are isolated, so $M'=\{B\neq0\}$ is open and dense. The level hypersurfaces of $\sigma$ in $M'$ are totally umbilical, with shape operator
		$$
		A=\frac{kf+C}{\|df\|}\,\mathbb{I},
		$$
		with respect to the unit normal $N=B/\|B\|$, and they have constant mean curvature.
		\item[(d)] Every point of $M'=\{B\neq 0\}$ has a neighbourhood
		isometric to a warped product
		$$
		\left(I\times\Sigma,\,
		dt^2+\left(\frac{f'(t)}{f'(t_0)}\right)^2g_\Sigma\right),
		$$
		where $\Sigma$ is the level hypersurface of $f$ through a fixed
		base point, equipped with its induced metric $g_\Sigma$, and
		$t$ is the signed distance from $\Sigma$ along the unit normal
		$$
		N=\frac{\mathrm{grad}_g f}{\|df\|}.
		$$
		The orientation of $t$ is chosen so that
		$$
		f'(t)=\frac{df}{dt}=\|df\|>0.
		$$
		Here $t_0$ is the value of the distance coordinate on the initial hypersurface $\Sigma$. The representation is local,
		on a neighbourhood where $df\neq 0$ and the normal flow of $N$ provides product coordinates. Along each normal geodesic, $f=f(t)$ satisfies
		$$
		f''(t)=-(kf(t)+C),\;\; f'(t)^2+kf(t)^2+2Cf(t)+1=0.
		$$		
	\end{enumerate}
\end{corollary}
\begin{proof}
	(a) follows from \eqref{eq:concircular} and \eqref{eq:nablaB}, since $\gamma+\|B\|^2=-k-C\alpha$, and (b) follows from \eqref{eq:Nk}. (c) At a critical point $p$ of $f$, \eqref{eq:firstint} gives $kf(p)^2+2Cf(p)+1=0$. If also $kf(p)+C=0$, then $kf(p)^2=-Cf(p)$, hence $Cf(p)=-1$ and $kf(p)^2=1$, so $C^2=k$, contradicting Corollary~\ref{cor:range}. Thus $\mathrm{Hess}_gf(p)=-(kf(p)+C)g$ is nondegenerate, and $p$ is an isolated critical point. On $M'$, let $N=\mathrm{grad}_gf/\|df\|=B/\|B\|$. By \eqref{eq:firstint}, $\|df\|$ is a function of $f$, hence constant on each level set. So for $X$ tangent to a level set, $X(\|df\|)=0$ and
	$
	\nabla_XN=\frac{1}{\|df\|}\nabla_X\mathrm{grad}_gf=-\frac{kf+C}{\|df\|}\,X .
	$
	Thus $A=-\nabla N=\frac{kf+C}{\|df\|}\mathbb{I}$, which is constant on each level set.
	
	(d)	On $M'=\{df\neq0\}$, the unit normal $N=\mathrm{grad}f/\|df\|$ is geodesic since
	$
	\mathrm{Hess}f=-(kf+C)g
	$
	implies $\nabla_N\mathrm{grad}f=-(kf+C)N$ and $\|df\|$ is a function of $f$ by the first-integral identity. Thus the integral curves of $N$ are unit-speed geodesics and provide local coordinates $(t,p)$ near any point of $M'$. For $X$ tangent to a level set,
	$
	\nabla_XN=-\frac{kf+C}{\|df\|}X.
	$
	Consequently the level hypersurfaces are totally umbilical, and their induced metrics evolve by a scalar factor along the normal flow. Writing $f=f(t)$, we have $f'=\|df\|$ and
	$
	f''=-(kf+C).
	$
	Integration of the metric evolution gives
	$$
	g=dt^2+\left(\frac{f'(t)}{f'(t_0)}\right)^2g_\Sigma.
	$$
	The first integral is the identity displayed in Theorem \ref{thm:concircular}. 
\end{proof}

\subsection{The complete case}

We denote by $\mathbb S^{2n+1}(k)$, $k>0$, the round sphere of curvature $k$ (and radius $1/\sqrt k$), and by $\mathbb H^{2n+1}(k)$, $k<0$, the hyperbolic space of curvature $k$.
\begin{theorem}\label{thm:complete}
	Let $M$ be a connected l.c.\ Sasakian $N(k)$-manifold with $\omega\neq 0$ whose metric $g$ is complete. Then $(M,g)$ is isometric to the simply connected space form of constant curvature $k$:
	\begin{enumerate}
		\item[(i)] the round sphere $\mathbb S^{2n+1}(k)$ if $k>0$;
		\item[(ii)] the Euclidean space $\mathbb{R}^{2n+1}$ if $k=0$;
		\item[(iii)] the hyperbolic space $\mathbb H^{2n+1}(k)$ if $k<0$.
	\end{enumerate}
\end{theorem}
\begin{proof}
By Theorem~\ref{thm:concircular}, the nonconstant function $f$ satisfies $\mathrm{Hess}_{g} f=-(kf+C)g$ and $\|df\|^2+kf^2+2Cf+1=0$. 
	
	(i) If $k>0$, then $F=f+C/k$ is nonconstant and satisfies $\mathrm{Hess}_gF=-kF\,g$. By Obata's theorem \cite{obata1962}, in its complete connected form, implies that $(M,g)$ is isometric to the round sphere of sectional curvature, i.e., $\mathbb S^{2n+1}(k)$.
	
	(ii) If $k=0$, then $\mathrm{Hess}_gf=-C\,g$. The first integral gives $\|df\|^2+2Cf+1=0$, so $C\neq0$, constant. By Tashiro's theorem \cite{tashiro}, $(M,g)$ is isometric to Euclidean space.
	
	(iii) If $k<0$, put $\kappa=\sqrt{-k}$ and $F=f+C/k$. Then $\mathrm{Hess}_gF=\kappa^2F\,g$ and, by \eqref{eq:firstint},
	$
	\|dF\|^2-\kappa^2F^2=-1-C^2/\kappa^2<0 .
	$
	In particular $F$ vanishes nowhere; let $\varepsilon=\pm1$ be its sign. Let $p$ be a regular point of $F$ and $\gamma$ the unit speed geodesic with $\gamma'(0)=-\varepsilon\,\mathrm{grad}\,F/\|dF\|$ at $p$. Since $\mathrm{grad}\,F$ is concircular, its integral curves are pregeodesics, so $\gamma'$ remains equal to $-\varepsilon\,\mathrm{grad}\,F/\|dF\|$ as long as $\mathrm{grad}\,F\neq0$ along $\gamma$. The function $F\circ\gamma$ satisfies $(F\circ\gamma)''=\kappa^2F\circ\gamma$, so $(F\circ\gamma)'^2-\kappa^2(F\circ\gamma)^2$ is constant along $\gamma$, equal to its negative value $\|dF\|^2(p)-\kappa^2F(p)^2$ at $t=0$. Hence $F\circ\gamma=A\cosh(\kappa(t-t_{0}))$ with $\varepsilon A>0$, and $t_{0}>0$ because $|F\circ\gamma|$ is decreasing at $t=0$. By completeness $\gamma$ is defined on $[0,t_{0}]$, and $\mathrm{grad}\,F$ must vanish at some point of $\gamma([0,t_{0}])$: otherwise $(F\circ\gamma)'(t)=-\varepsilon\|dF\|(\gamma(t))\neq0$ on $[0,t_{0}]$, contradicting $(F\circ\gamma)'(t_{0})=0$. Thus $F$ has a critical point. A complete Riemannian manifold carrying a nonconstant solution of $\mathrm{Hess}\,F=\kappa^2F\,g$ with a critical point is isometric to the hyperbolic space of curvature $-\kappa^2$ \cite{tashiro,kanai}. Hence $(M,g)$ is isometric to $\mathbb H^{2n+1}(k)$.
\end{proof}
Theorem~\ref{thm:complete} shows that the $k$-nullity condition, which a priori constrains only the curvature in the direction of $\xi$, determines the entire Riemannian geometry of a complete l.c.\ Sasakian manifold with nonzero Lee form. All three cases occur (Remark~\ref{rem:spaceforms}). Without completeness, the condition still allows manifolds that are not of constant curvature in every dimension $2n+1\ge5$ (Example~\ref{ex:join}); the three-dimensional case is treated in \cite{partII}.

\subsection{Examples}\label{sec:examples}

The examples below come from the following construction. Recall that a function $v$ on a Riemannian manifold $(M_{0},g_{0})$ satisfies \emph{Obata's equation} if $\mathrm{Hess}_{g_{0}}v=-v\,g_{0}$ \cite{obata1962}.
\begin{proposition}\label{prop:construction}
	Let $(M_{0},\phi_{0},\xi_{0},\eta_{0},g_{0})$ be a connected Sasakian manifold, and let $v$ be a non-constant smooth solution of Obata's equation on $M_{0}$. Then $c^2:=v^2+\|dv\|^2_{g_{0}}$ is a positive constant. Let $C\in\mathbb{R}$, let $U$ be a connected component of $\{v>C\}$, and put $u=v-C>0$ and $\sigma=-\ln u$. Then the l.c.\ Sasakian structure of Example~\ref{ex:lcs} on $U$,
	$$
	\phi=\phi_{0},\quad \xi=u\,\xi_{0},\quad \eta=u^{-1}\eta_{0},\quad g=u^{-2}g_{0},
	$$
	has the following properties.
	\begin{enumerate}
		\item[(i)] $\omega=-du/u\neq 0$, and $(U,\phi,\xi,\eta,g)$ is an l.c. Sasakian $N(k)$-manifold with $k=C^2-c^2$. The function $f=1/u$ satisfies \eqref{eq:concircular} with the constant $C$.
		\item[(ii)] $R_g=k\,G_g+u^{-2}(R_{g_{0}}-G_{g_{0}})$. Equivalently, $K_g(\pi)=k+u^2\big(K_{g_{0}}(\pi)-1\big)$ for every plane $\pi$, where $K_g$ and $K_{g_{0}}$ are the sectional curvatures. In particular, $(U,g)$ has constant curvature if and only if $(U,g_{0})$ has constant curvature $1$.
	\end{enumerate}
\end{proposition}
\begin{proof}
	Since $d\big(v^2+\|dv\|^2_{g_{0}}\big)=2v\,dv+2\mathrm{Hess}_{g_{0}}v(\mathrm{grad}_{g_{0}}v,\cdot)=0$, $c^2$ is constant, and $c>0$ because $v$ is not constant. Along every $g_{0}$-geodesic $\tau$ we have $(v\circ\tau)''=-v\circ\tau$; hence, if $v$ and $dv$ vanished at a point, $v$ would vanish identically. If $v$ were constant on an open set, Obata's equation would force $v=0$ there, so $v\equiv0$. Thus $v$ is not constant on any open set, and $\omega=d\sigma=-du/u\neq 0$ on $U$. Here $f=\exp(\sigma)=1/u$ and $\alpha=u$.
	
	(i) By \eqref{eq:conn}, for every function $h$,
	$$
	\mathrm{Hess}_gh(X,Y)=\mathrm{Hess}_{g_{0}}h(X,Y)-\omega(X)\,dh(Y)-\omega(Y)\,dh(X)+g(X,Y)\,dh(B).
	$$
	For $h=f$ we have $df=f\omega$ and $df(B)=f\|B\|^2$, and, since $du=-u\,\omega$, $\mathrm{Hess}_{g_{0}}f=-u^{-2}\mathrm{Hess}_{g_{0}}u+2f\,\omega\otimes\omega$. Hence
	$$
	\mathrm{Hess}_gf=-u^{-2}\mathrm{Hess}_{g_{0}}u+f\|B\|^2g .
	$$
	Now $\mathrm{Hess}_{g_{0}}u=\mathrm{Hess}_{g_{0}}v=-(u+C)g_{0}=-(u+C)u^2g$ and $\|B\|^2=u^2\|\omega\|^2_{g_{0}}=\|du\|^2_{g_{0}}=c^2-(u+C)^2$. Therefore
	$$
	\mathrm{Hess}_gf=\Big\{(u+C)+\frac{c^2-(u+C)^2}{u}\Big\}g=\frac{c^2-C^2-Cu}{u}\,g=-(kf+C)\,g ,
	$$
	with $k=C^2-c^2$. By Theorem~\ref{thm:concircular}, applied to the connected manifold $U$, the structure is an l.c.\ Sasakian $N(k)$-manifold.
	
	(ii) The curvature of $g=u^{-2}g_{0}$ is given by the transformation law under conformal changes \cite[Theorem~1.159]{besse}, which in our conventions reads
	$$
	R_g=u^{-2}\Big(R_{g_{0}}+u^{-1}\,\mathrm{Hess}_{g_{0}}u\KN g_{0}-\tfrac12u^{-2}\|du\|^2_{g_{0}}\,g_{0}\KN g_{0}\Big),
	$$
	where $\KN$ is the Kulkarni--Nomizu product, so that $G_h=\frac12 h\KN h$. Since $\mathrm{Hess}_{g_{0}}u=-(u+C)g_{0}$, $\|du\|^2_{g_{0}}=c^2-(u+C)^2$ and $G_{g_{0}}=u^4G_g$, this gives
	$$
	R_g=u^{-2}(R_{g_{0}}-G_{g_{0}})+\big\{u^2-2u(u+C)-c^2+(u+C)^2\big\}G_g=u^{-2}(R_{g_{0}}-G_{g_{0}})+k\,G_g .
	$$
	The formula for sectional curvatures follows from $G_{g_{0}}=u^4G_g$. If $K_{g_{0}}=1$, then $K_g= k$. Conversely, if $K_g= c'$ is constant, then $K_{g_{0}}=1+(c'-k)u^{-2}$ is pointwise constant, hence constant by Schur's lemma, since $\dim U\ge3$, as $u$ is not constant, $c'=k$ and $K_{g_{0}}=1$.
\end{proof}

Conversely, every l.c.\ Sasakian $N(k)$-manifold with $\omega\neq 0$ arises in this way from its underlying Sasakian structure, with $v=\exp(-\sigma)+C$; this equivalence and its consequences are studied in \cite{partII}. 
\begin{remark}\label{rem:spaceforms}
	\emph{Let $(\mathbb S^{2n+1},\phi_{0},\xi_{0},\eta_{0},g_{0})$ be the standard Sasakian unit sphere in $\mathbb C^{n+1}=\mathbb{R}^{2n+2}$, and let $v(x)=\langle x,a\rangle$ with $a\in\mathbb{R}^{2n+2}\setminus\{0\}$, so that $\mathrm{Hess}_{g_{0}}v=-v\,g_{0}$ and $c=|a|$. For $C<|a|$, Proposition~\ref{prop:construction} applies on the connected set $U=\{v>C\}$, and since $R_{g_{0}}=G_{g_{0}}$, the metric $g=(v-C)^{-2}g_{0}$ has constant curvature $k=C^2-|a|^2$. This metric is complete, and $(U,g)$ is a round sphere if $C<-|a|$, a Euclidean space if $C=-|a|$ (by stereographic projection from $-a/|a|$), and a hyperbolic space if $-|a|<C<|a|$; in the last case $U$ is a spherical cap and $v-C$ vanishes to first order on $\partial U$, since $\|dv\|^2_{g_{0}}=|a|^2-C^2>0$ there. Thus all three cases of Theorem~\ref{thm:complete} occur.}
\end{remark}
\begin{example}\label{ex:join}
	\emph{Let $n\ge2$ and let $(N,\phi_N,\xi_N,\eta_N,g_N)$ be $\mathbb{R}^{2n-1}$ with the standard Sasakian structure of Example~\ref{ex:R2n+1} (with $n-1$ in place of $n$), which has constant $\phi$-sectional curvature $-3$ \cite{bl}. On
		$
		L=(0,\pi/2)\times\mathbb{R}\times N,
		$
		with coordinates $(t,s)$ on the first two factors, consider
		\begin{align*}
			g_{0}&=dt^2+\cos^2t\,ds^2+\sin^2t\,g_N,\;\; \xi_{0}=\frac{\partial}{\partial s}+\xi_N,\;\; \eta_{0}=\cos^2t\,ds+\sin^2t\,\eta_N,\\
			\phi_{0}&=\phi_N+\Big(\tan t\,\frac{\partial}{\partial s}-\cot t\,\xi_N\Big)\otimes dt+\sin t\cos t\,\frac{\partial}{\partial t}\otimes(\eta_N-ds),
		\end{align*}
		where $\phi_N,\xi_N,\eta_N,g_N$ are pulled back to $L$, so that $\phi_N(\partial/\partial t)=\phi_N(\partial/\partial s)=0$. One checks directly from \eqref{equa1}, \eqref{equation2} and \eqref{eq:sasakian} that $(\phi_{0},\xi_{0},\eta_{0},g_{0})$ is a Sasakian structure on $L$. Geometrically, with $w=\cos t\,\exp(\mathrm{i}s)$ and $r=\sin t$, the manifold $(L,g_{0})$ is locally isometric to the link $\{|w|^2+r^2=1\}$ of the K\"ahler cone $\mathbb C\times C(N)$, where $C(N)=(\mathbb{R}_{+}\times N,dr^2+r^2g_N)$ is the K\"ahler cone of $N$ \cite{boyer-galicki}. The function
		$$
		v=\cos t\cos s,
		$$
		which is the restriction of $\mathrm{Re}\,w$, satisfies $\mathrm{Hess}_{g_{0}}v=-v\,g_{0}$ and $v^2+\|dv\|^2_{g_{0}}=1$.\\
		Let $C\in\mathbb{R}$, and let $U$ be a connected component of $\{v>C\}$; for $C<-1$, $U=L$. By Proposition~\ref{prop:construction}, $g=(v-C)^{-2}g_{0}$ defines an l.c.\ Sasakian $N(k)$-structure on $U$ with $k=C^2-1$. The metric $g_{0}$ is a multiply warped product, and a direct computation, as for warped products, shows that for vectors $X,Y$ tangent to $N$,
		$$
		K_{g_{0}}(X,Y)=\frac{K_N(X,Y)-\cos^2t}{\sin^2t}.
		$$
		The plane spanned by $X_1=\partial/\partial x_1+y_1\,\partial/\partial z$ and $X_2=\partial/\partial y_1$ is $\phi_N$-invariant, since $\phi_NX_1=-X_2$, so $K_N(X_1,X_2)=-3$ and $K_{g_{0}}(X_1,X_2)=-(3+\cos^2t)/\sin^2t$. By Proposition~\ref{prop:construction}(ii),
		$$
		K_g(X_1,X_2)= k-\frac{4(\cos t\cos s-C)^2}{\sin^2t}.
		$$
		Hence $(U,g)$ is not of constant curvature, although every plane containing $\xi$ has sectional curvature $k$. By Theorem~\ref{thm:complete}, $g$ is not complete. Since $n\ge2$ is arbitrary, such examples exist in every dimension $2n+1\ge5$, and for every sign of $k$.
	}
\end{example}

\section{Concluding remarks}\label{sec:conclusion}

The results above show that, after a conformal change of a Sasakian structure, the $k$-nullity condition is no longer a mild curvature assumption. It is equivalent to the fact that $g$ and the Sasakian metric $g_{0}$ have the same geodesic circles, it forces the inverse conformal factor to be a concircular function, it determines the range of the conformal factor, it forces $k>0$ on compact manifolds, and it determines the Riemannian geometry completely when the metric is complete. Several questions remain.
\begin{enumerate}
	\item[(1)] Proposition~\ref{prop:construction} shows that solutions of Obata's equation on Sasakian manifolds produce l.c.\ Sasakian $N(k)$-manifolds. The converse also holds. It leads to a description of all local models through K\"ahler cones and to a totally geodesic foliation spanned by $\xi$, $B$ and $\phi B$; this is the subject of the companion paper \cite{partII}.
	\item[(2)] By Remark~\ref{rem:convention}, a $(\kappa,\mu)$-type condition with a non-trivial $\mu$-term requires an operator that is symmetric on $\ker\eta$. Natural candidates are built from the symmetric part of $\nabla B$. It would be interesting to know which of the resulting conditions admit non-trivial examples.
	\item[(3)] Since l.c.\ Sasakian structures always have exact Lee forms, a genuinely local theory requires different local models or different deformations, for instance $\mathcal D$-homothetic deformations of Sasakian structures, which are not conformal (Remark~\ref{rem:global}).
	\item[(4)] Further curvature conditions, such as $\eta$-Einstein or Ricci soliton conditions on l.c.\ Sasakian manifolds, should be accessible through the identities of Section~\ref{sec:basic} and the concircular structure of Section~\ref{sec:consequences}.
\end{enumerate}

\section*{Acknowledgments} 
This work is partially supported by the National Research Foundation of South Africa.

%\section*{Declarations} 

%\noindent\textbf{Data availability.} Data sharing is not applicable to this article, as no datasets were generated or analysed.

%\smallskip
%\noindent\textbf{Conflict of interest.} The authors declare that they have no conflict of interest.

\end{document}